\documentclass{primus}
\usepackage{graphicx}
\usepackage{amsmath}
\usepackage{amsfonts}
\usepackage{amssymb}
\usepackage{color}
\usepackage{amsthm}
\usepackage{float}
\usepackage[utf8]{inputenc}
\usepackage[english]{babel}

\newcommand{\AmSLaTeX}{$\cal A$\kern-.1667em\lower.5ex\hbox{$\cal
M$}\kern-.125em $\cal S$-\LaTeX}

\newtheorem{theorem}{Theorem}

\begin{document}

\title {AREA INSIDE THE CIRCLE:\\ INTUITIVE AND RIGOROUS PROOFS}

\author {V. Siadat \\
Department of Mathematics\\
Richard J. Daley College\\
Chicago, IL 60652,  USA\\
vsiadat@ccc.edu\\
(773) 838-7658\\}

\markboth{V. Siadat}{AREA INSIDE THE CIRCLE}


\makePtitlepage

\begin{abstract}
 In this article I conduct a short review of the proofs of the area inside a circle. These include intuitive as well as rigorous analytic proofs. This discussion is important not just from mathematical view point but also because pedagogically the calculus books still today use circular reasoning to prove the area inside a circle (also that of an ellipse) on  this important historical topic, first illustrated by Archimedes. I offer an innovative approach, through introduction of a theorem, which will lead to proving the area inside a circle avoiding circular argumentation.
\end{abstract}

\listkeywords

\section*{ACKNOWLEDGMENTS}

I wish to thank Professor Cyrill Oseledets of Richard J. Daley College for his review of the initial manuscript and making helpful suggestions to improve it.

\newpage

\section{INTRODUCTION}
  Why area inside a circle again? Why we shouldn't confound the notions of
intuition and rigor? Do calculus books, even today, still resort to circular
reasoning? This paper is an attempt to elucidate these questions by walking
the reader through the path of intuitive to solid analytical reasoning,
pointing out the gaps that often occur, on the proof of this ancient and well
known problem, first illustrated by Archimedes. The motivation behind writing
of this piece was to engage the reader in further thinking about mathematical
proofs and the level of rigor at which they are presented.

In the following we present a brief review of the proofs of area inside a
circle. A typical rigorous proof requires knowledge of integral calculus, see
for example \cite{Larson}. But even in these proofs presented by calculus books, see for
example \cite{Briggs}, the authors resort to circular reasoning. To prove the area
inside a circle, they set up the integral  $\int_{0}^{1}\sqrt{1-x^{2}}dx$ 
followed by trigonometric substitution which requires
knowing that the derivative of $\sin\theta$ is $\cos\theta.$ But this latter
fact requires proving that $\lim_{\theta\rightarrow0}\dfrac{\sin\theta}%
{\theta}=1.$ For this proof, they resort to a geometric argument, bounding the
area of a sector of a unit circle between the areas of two triangles and
showing that $\sin\theta<\theta<\tan\theta$. They then apply the Squeeze
Theorem. But for computation of the sector's area, they resort to a standard
formula, $A=\frac{1}{2}\theta,$ which is based on knowing the area of a
circle. So, they prove the area by assuming the area. This is obviously
circular argumentation! For an excellent critique of this method see \cite{Krantz}.
There are also a number of intuitive proofs intended to provide an insight to
the derivation of the area with just the knowledge of geometry and limits. In
this short piece we begin by proving a preliminary result showing that
$\lim_{\theta\rightarrow0}\dfrac{\sin\theta}{\theta}=1,$ without, a priori,
assuming the area of a sector. This limit is central to the proof of the
derivatives of trigonometric functions. We note that aside from the
aforementioned limit, the function $\dfrac{\sin\theta}{\theta}$ itself plays
an important role not only in mathematics but in other fields of science such
as physics and engineering.
\section{PROOFS}
Consider a circle of radius 1, centered at the origin, as shown in Fig. 1; see
\cite{Siadat}.

\begin{figure}[H]
\centering
\includegraphics[scale=0.6]{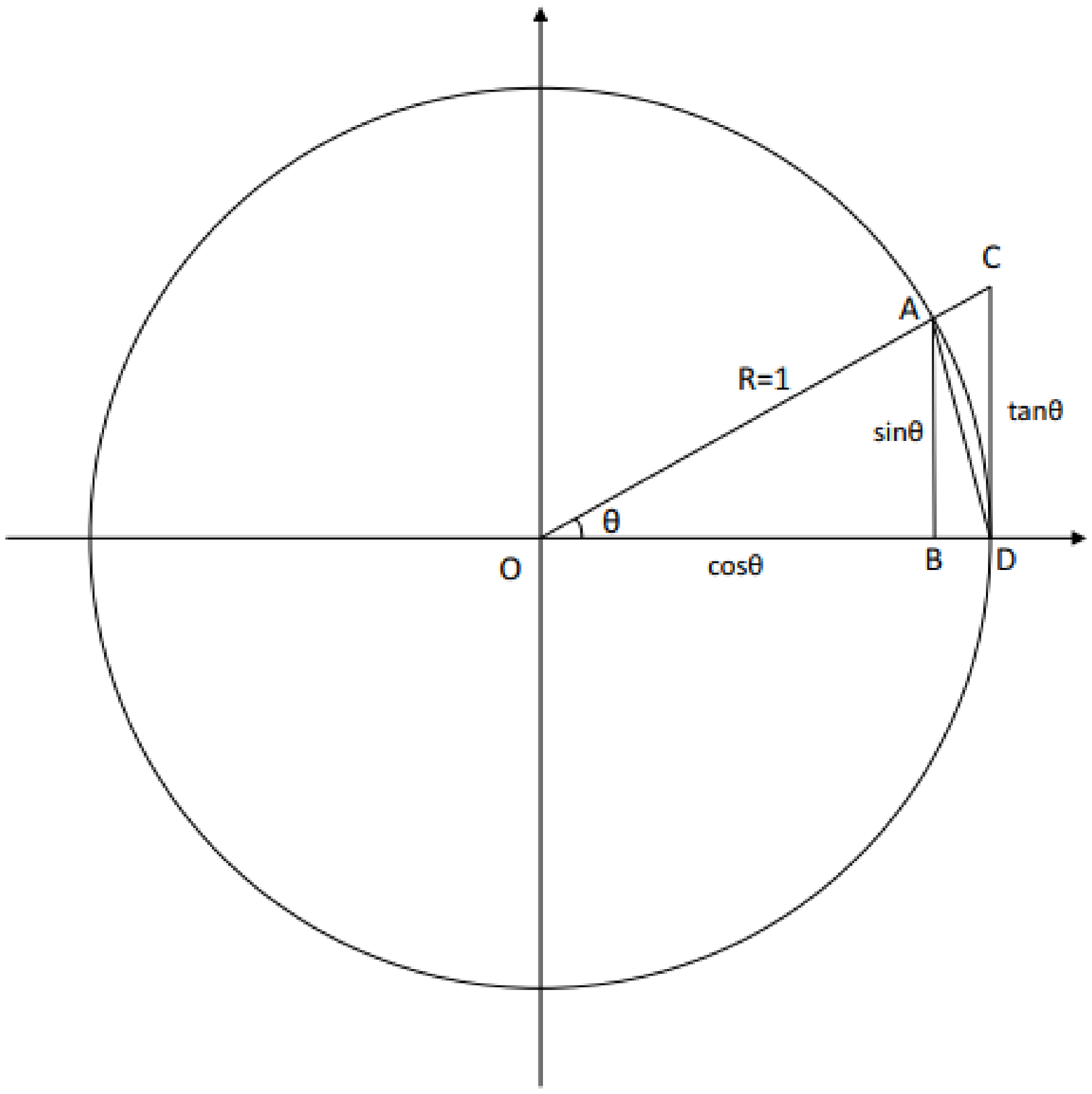}
\caption{ \label{figure:1}}
\end{figure}

\begin{theorem}
Let $0<\theta<\dfrac{\pi}{2}$ be an angle measured in radians.
Then, 
\begin{equation*}
\lim_{\theta\rightarrow0}\dfrac{\sin\theta}{\theta}=1.
\end{equation*}
\end{theorem}

\begin{proof} 
 Since the magnitude of $\theta$   equals the length of the arc it subtends and since $\sin\theta<\overline{AD},$ we have that $\sin\theta<\theta,$ or $1<\dfrac{\theta}{\sin\theta}.$ This establishes a lower bound for $\dfrac{\theta}{\sin\theta}.$ To show the upper bound, observe that by the triangle inequality, $\theta<\sin\theta+\overline{BD}.$ This can be established by the standard method of estimating an arc length of a rectifiable curve by the linear approximation of the lengths of the chords it subtends through partitioning. The result follows by applying the triangle inequality in each partition; see \cite{Krantz}. Noting that 
 $\sin\theta<\tan\theta$
 we get $\theta<\tan\theta+\overline{BD}$
 But $\overline{BD}=1-\cos\theta.$ So, 
 \begin{equation*}
 \dfrac{\theta}{\sin\theta}<\dfrac{1}{\cos\theta}+\dfrac{1-\cos\theta} {\sin\theta}
 \end{equation*}
 and 
 \begin{equation*}
 \dfrac{1-\cos\theta}{\sin\theta}=\sqrt{\dfrac{(1-\cos \theta)^{2}}{1-\cos^{2}\theta}}=\sqrt{\dfrac{1-\cos\theta}{1+\cos\theta}}
 \end{equation*}
 Combining this result with the previous lower bound gives, 
 \begin{equation*}
 1<\dfrac{\theta}{\sin\theta}<\dfrac{1}{\cos\theta}+\sqrt{\dfrac{1-\cos\theta}{1+\cos\theta}}
 \end{equation*}
 Letting $\theta\rightarrow0$ in the last expression completes the upper bound, resulting in $1 \leq \dfrac{\theta}{\sin\theta} \leq 1.$ \ Finally, applying the Squeeze Theorem we get,

\begin{equation} \label{limit} 
\lim_{\theta\rightarrow 0}\dfrac{\sin\theta}{\theta}=1
\end{equation}
\end{proof}

An interesting question related to the foregoing bounding of the angle
$\theta$ is that if we define the derivatives of trigonometric functions of
$\theta$ analytically (i.e., by infinite series or complex numbers or solutions
of differential equations), can we arrive at the bounding of the angle? The
next theorem follows.

\begin{theorem} If $(\sin\theta)^{^{\prime}}=\cos\theta,$ and $(\tan
\theta)^{^{\prime}}=\sec^{2}\theta,$ then

$\sin\theta<\theta<\tan\theta.$

\end{theorem}
\begin{proof}

Let $f(\theta)=\sin{(\theta)}-\theta$ Than $f(0)=0$ and $f^{^{\prime}}(\theta)=\cos\theta-1\leq0,$ giving
$f^{^{\prime}}(\theta)<0$ for $\theta>0.$

Hence,

$f(\theta)=%
%TCIMACRO{\dint _{0}^{\theta}}%
%BeginExpansion
{\displaystyle\int_{0}^{\theta}}
%EndExpansion
f^{^{\prime}}(t)dt<0.$

Therefore, $\sin\theta<\theta.$ This gives a lower bound for $\theta$.

To find an upper bound for $\theta,$ let $g(\theta)=\tan\theta-\theta.$

$g(0)=0$ and $g^{^{\prime}}(\theta)=\sec^{2}\theta-1\geq0,$ giving
$g^{^{\prime}}(\theta)>0$ for $\theta>0.$

Hence,

$g(\theta)=%
%TCIMACRO{\dint _{0}^{\theta}}%
%BeginExpansion
{\displaystyle\int_{0}^{\theta}}
%EndExpansion
g^{^{\prime}}(t)dt>0.$

Therefore, $\theta<\tan\theta.$ This gives an upper bound for $\theta$.
Combining the above results we get $\sin\theta \leq \theta \leq \tan\theta.$

\end{proof}

\begin{theorem}

 Area inside a circle of radius $R$ is $\pi R^{2}.$

\end{theorem}

\begin{proof}

Consider a circle of radius $R$ centered at the origin in Fig.
2. 
\begin{figure}[H]
\centering
\includegraphics[scale=0.6]{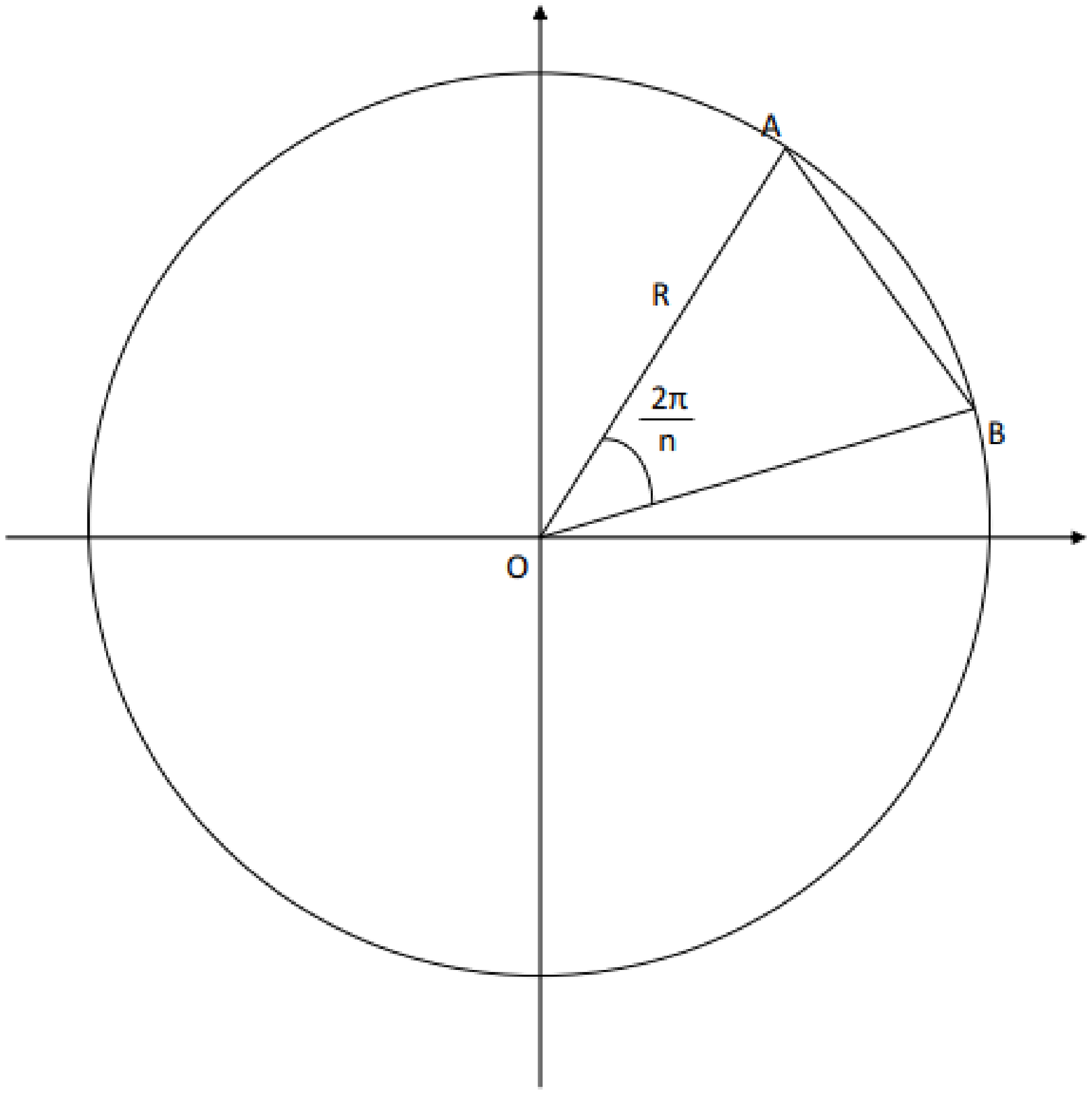}
\caption{ \label{figure:sine}}
\end{figure}

Partition the circle into $n$ equal slices and consider a slice with
central angle $\frac{2\pi}{n}$ radians. We know that the area of a triangle is
one-half times the product of two of its sides times the sine of the angle
between the two sides.\ So, the area of the triangle subtended by the central
angle $\frac{2\pi}{n}$ becomes $A=\frac{1}{2}R^{2}\sin(\frac{2\pi}{n}).$
Because there are $n$ inscribed triangles in the circle, the total area of all
these triangles would be $A_{\text{total}}=n\frac{1}{2}R^{2}\sin(\frac{2\pi
}{n}).$ As we increase the number of slices by increasing $n,$ the sum of the
areas of the inscribed triangles get closer to the area of the circle. To get
the area of the circle, we need to find the limit of $A_{\text{total}}$, as
$n\rightarrow\infty.$ So, using (\ref{limit}), and since $\frac{2\pi}{n}%
\rightarrow0,$ as $n\rightarrow\infty,$ the area of the circle becomes:

$A_{\text{circle}}=\lim_{n\rightarrow\infty}$ $n\frac{1}{2}R^{2}\sin
(\frac{2\pi}{n})=\pi R^{2}(\lim_{n\rightarrow\infty}\dfrac{\sin\frac{2\pi}{n}%
}{\frac{2\pi}{n}})=\pi R^{2}.$

Hence, $A_{\text{circle}}=\pi R^{2}.$
\end{proof}
An intuitive and interesting method of proving the area inside a circle which
requires area stretching \footnote{Area stretching is a result from geometry
stating that if we stretch a region in the coordinate plane vertically by a
factor of $k>0$ and horizontally by a factor of $l>0,$ then its area will
stretch by the factor $kl.$} and mapping from an annulus to a trapezoid is
discussed in \cite{Axler}. The subtle point in this method, as expressed in \cite{Axler}, is
that it assumes as evident the area preservation from a circular to a simply
connected region. For a discussion of transformation of different regions
using complex variable method, see \cite{Brown}.

There are many other intuitive approaches also, some of which involve slicing
or opening up a circle. Below we offer a simple intuitive proof which is not
based on area stretching, but assumes area preservation under mappings.
Consider two concentric circles with radii $r$ and $R$ and corresponding areas
$A_{r\text{ }}$and $A_{R}$. Cut the annulus open in the shape of a right angle
trapezoid $ABCD$ as in Fig. 3.

\begin{figure}[H]
\centering
\includegraphics[scale=0.80]{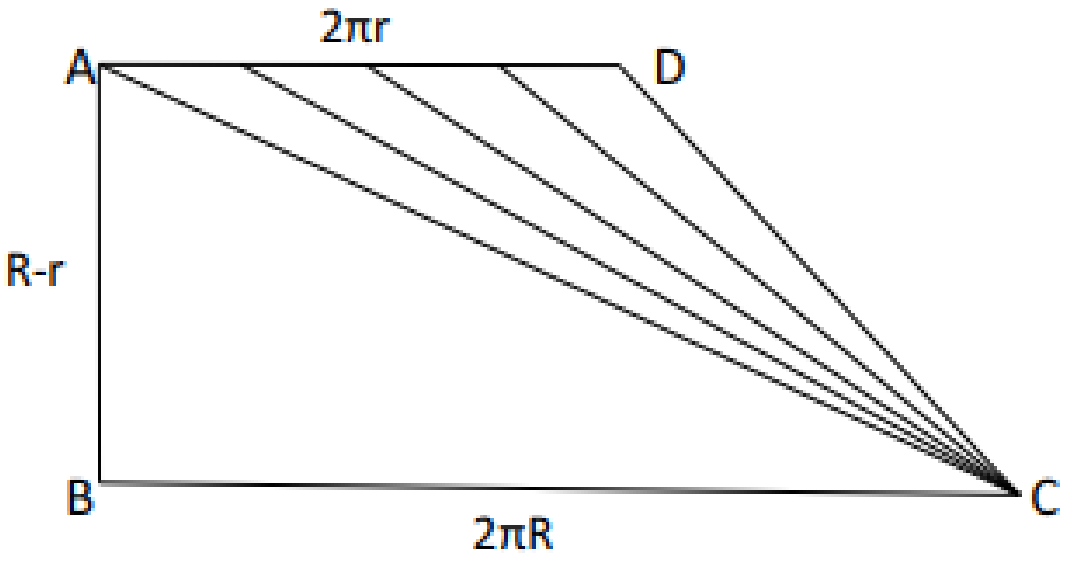}
\caption{ \label{figure:3}}
\end{figure}

 We can see that the area of the annulus equals the area of the
trapezoid. So, $A_{\text{annulus}}=A_{\text{trapezoid}}$, or \begin{equation*}
A_{R}-A_{r}=\frac{1}{2}(2\pi R+2\pi r)(R-r)=\pi R^{2}-\pi r^{2}.
\end{equation*}
We can choose
$r>0$ as small as we please and so, in particular, if we let $r$ approach $0,$
the area of the inner circle approaches $0$ and we get, $A_{R}-0=\pi R^{2}-0.$
Thus $A_{R}=\pi R^{2}.$
Note that shrinking $r$ to $0,$ shrinks the trapezoid to the right triangle
$ABC,$ whose area is $\frac{1}{2}(2\pi R)(R-0)=\pi R^{2}=A_{R}.$
In the following we present an analytic proof of the area inside a circle
using area stretching, which does not assume area preserving mapping of regions.
\begin{theorem}
Area inside a circle of radius $r$ is $\pi r^{2}.$
\end{theorem}
\begin{proof}
 Consider a circle of radius $r$ centered at the origin and partition it into $n$ equal sectors, each having central angle $\dfrac{2\pi}{n},$ and the corresponding arc length $\dfrac{2\pi}{n}r.$ Assume the area of a sector is $c_{n}.$ If $\ $we stretch the radius $r$ by a factor of $k>1,$ we create a circle with radius $R=kr.$ So, the corresponding streched sector will have an arc length equal to $\dfrac{2\pi}{n}kr$ and the its area will be increased by a factor of $k^{2}$ to $k^{2}c_{n};$ see Fig. 4\ \ Now the area between the two sectors is $A_{\text{between sectors}}=k^{2}c_{n}-c_{n}=c_{n}[k^{2}-1]$ which is approximately equal to the area of the trapezoid $ABCD,$ in Fig. 4.
\begin{figure}[H]
\centering
\includegraphics[scale=0.6]{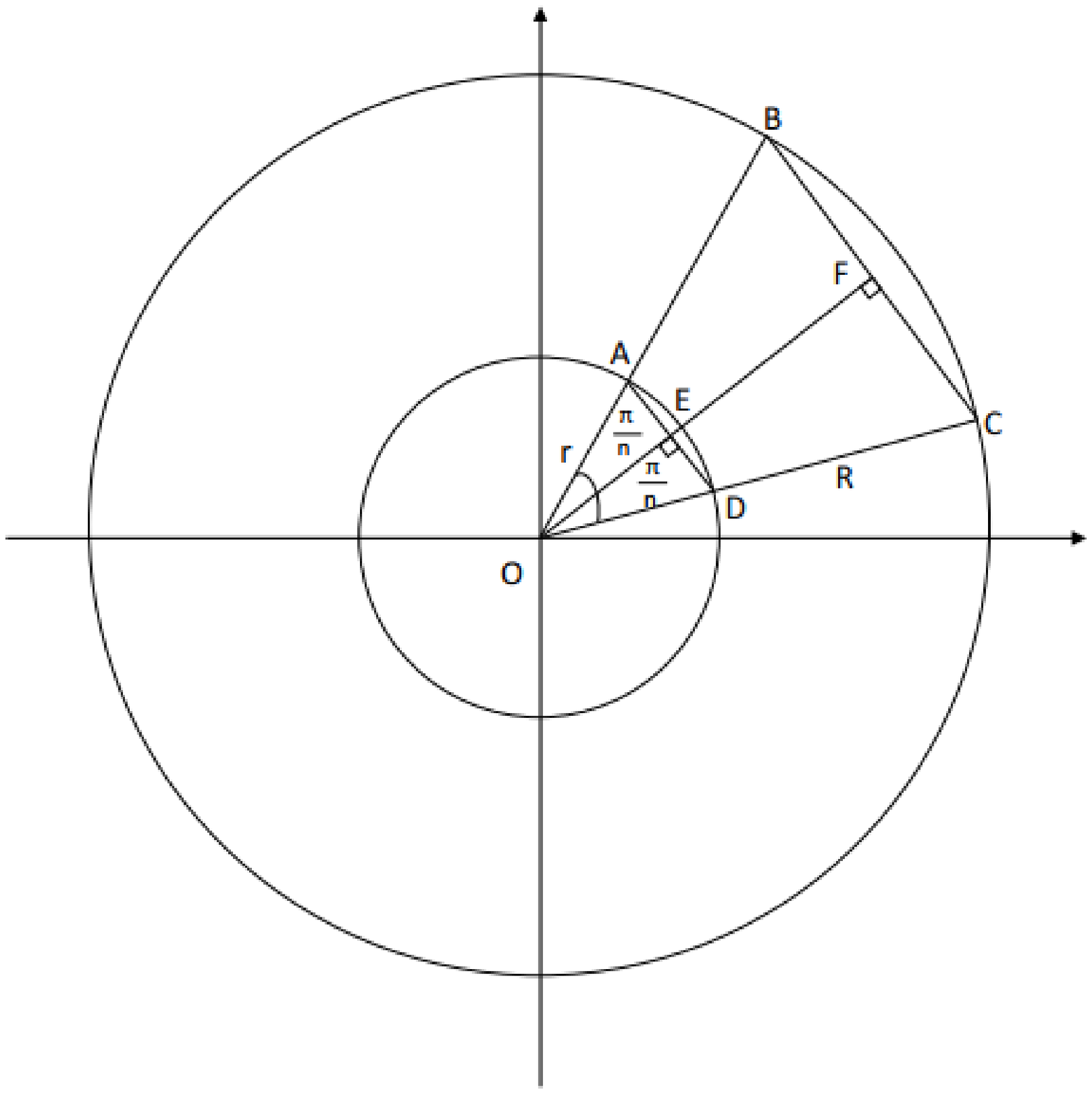}
\caption{ \label{figure:sine}}
\end{figure}

If we connect the center $O$ to the point $F$ which is the midpoint of
$\overline{BC},$ the triangles $\Delta OBF$ and $\Delta OCF$ become right
angle congruent triangles with right angles at the point $F.$ As a result,
central angles $\angle AOE$ and $\angle DOE$ will each equal $\dfrac{\pi}{n}.$
Obviously triangles $\Delta OAE$ and $\Delta ODE$ are also congruent having
right angles at the point $E.$ To calculate the area of the trapezoid, we note
that its larger base has length $\overline{BC}=2R\sin\dfrac{\pi}{n},$ and its
smaller base has length $\overline{AD}=2r\sin\dfrac{\pi}{n}.$ The height of
the trapezoid is: 
\begin{equation*}
\overline{EF}=R \cos\dfrac{\pi}{n}-r\cos\dfrac{\pi}{n}
\end{equation*}
Therefore the area of the trapezoid becomes:
\begin{eqnarray*}
A_{\text{trapezoid}}&=&\frac{1}{2}\left[2R\sin\dfrac{\pi}{n}+2r\sin\dfrac{\pi}%
{n}\right]\cdot\left[ R\cos\dfrac{\pi}{n}-r\cos\dfrac{\pi}{n}\right]\\
&=&\sin\dfrac{\pi}{n}\cos\dfrac{\pi}{n}[R^{2}-r^{2}]\notag\\
&=&\sin\dfrac{\pi}{n}\cos\dfrac{\pi}{n}\cdot r^{2}[k^{2}-1].\notag
\end{eqnarray*}
Setting $A_{\text{between sectors}}\approx A_{\text{trapezoid}}$, gives,
\begin{equation*}
\ c_{n}[k^{2}-1]\approx\sin\dfrac{\pi}{n}\cos\dfrac{\pi}{n}\cdot r^{2}%
[k^{2}-1],
\end{equation*} or $c_{n}\approx\sin\dfrac{\pi}{n}\cos\dfrac{\pi}{n}\cdot r^{2}.$ This approximation can be improved by increasing $n.$ Now, multiplying both
sides of the above by $n$ gives:
\begin{equation*}
nc_{n}\approx n\sin\dfrac{\pi}{n}\cos\dfrac{\pi}{n}\cdot r^{2}
\end{equation*} Since there are exactly $n$ identical sectors in the circle of radius $r,$ its
area becomes $c=nc_{n}.$ Therefore, $c$ $\approx n\sin\dfrac{\pi}{n}\cos\dfrac{\pi}{n}\cdot r^{2}.$ Now, taking the limit of both sides as 
$n\rightarrow\infty,$ and applying our earlier result (\ref{limit}) and the fact that $\cos\theta\rightarrow1,$ as $\theta\rightarrow0$,we get:

\begin{eqnarray*}
c&=&\lim_{n\rightarrow\infty}n\sin\dfrac{\pi}{n}\cos\dfrac{\pi}{n}\cdot r^{2}\\
&=&\lim_{n\rightarrow\infty}\pi\cdot\dfrac{\sin\frac{\pi}{n}}{\frac{\pi}
{n}}\cos\dfrac{\pi}{n}\cdot r^{2}\\
&=&\pi r^{2}(\lim_{n\rightarrow\infty}\dfrac{\sin\frac{\pi}{n}}{\frac{\pi}%
{n}})\cdot(\lim_{n\rightarrow\infty}\cos\dfrac{\pi}{n})\\
&=&\pi r^{2}
\end{eqnarray*}
Hence, $c=\pi r^{2}.$
\end{proof}

\pagebreak

\section*{BIOGRAPHICAL SKETCHES}

M. Vali Siadat is distinguished professor of mathematics at Richard J. Daley College. He holds two doctorates in mathematics, a Ph.D. in pure mathematics (harmonic analysis) and a D.A. in mathematics with concentration in mathematics education. Dr. Siadat has extensive publications in mathematics and mathematics education journals and has had numerous presentations at regional and national mathematics meetings. He is the recipient of the Carnegie Foundation for the Advancement of Teaching Illinois Professor of the Year Award in 2005 and the Mathematical Association of America's Deborah and Franklin Tepper Haimo Award in distinguished teaching of mathematics in 2009.

\end{document}